\documentclass[12pt,a4paper]{article}
\usepackage[utf8]{inputenc}
\usepackage{amsmath}
\usepackage{amsfonts}
\usepackage{amssymb}
\usepackage{makeidx}
\usepackage{graphicx}

\usepackage[left=2.25cm, right=2.25cm, top=2.25cm, bottom=2.25cm]{geometry}

\usepackage{moreverb}
\usepackage{bbm}

\newtheorem{d1}{Definition}
\newtheorem{d2}[d1]{Definition}

\newtheorem{p1}[d1]{Property}

\DeclareMathOperator\Sig{Sig}
\DeclareMathOperator\Tanh{Tanh}

\title{Abstract, keywords and references template}
\author{M. B\l asik$^{1}$  \\
	\small $^{1}$Institute of Mathematics, Czestochowa University of Technology \\
	Armii Krajowej 21, 42-201 Czestochowa, Poland \\
	marek.blasik@gmail.com
}
\date{} 

\title{Numerical method for the one phase 1D fractional Stefan problem 
	supported by an artificial neural network}

\begin{document}
\maketitle

\abstract{In this paper we present a numerical solution of a one-phase 1D fractional Stefan
	problem with Caputo derivative with respect to time variable. In the proposed approach,
	we use a front fixing method and the algorithm of numerical integration supported by 
	an artificial neural network. In the final part, we present some examples illustrating the 
	comparison of the new numerical scheme with its previous version and approximate analytical 
	solution.}

\section{Introduction}
Diffusion is a physical process that refers to the movement of molecules from a region of
high concentration to one of lower concentration. The phenomenon of classical diffusion is well 
described by the following relationship
\begin{equation}
\label{eq:diff}
<x^2(t)>\sim D t,
\end{equation}
where $<x^2(t)>$ is a mean squared displacement of the diffusing molecule in the course of time $t$,
$D$ is the diffusion coefficient. However, it turns out that in many phenomena occurring in nature,
the relation (\ref{eq:diff}) is not sufficient to describe them, which is confirmed by the results of 
experiments in \cite{Wee96,Sol93,Hum10,Kos05a,Kos05b}. 
A more general process of random motion of molecules is the so-called anomalous diffusion, which is 
characterized by the  non-linear relationship \cite{Met00,Met04}:
\begin{equation}
\label{eq:subdiff}
<x^2(t)>\sim D_{\alpha} t^{\alpha},
\end{equation}
where $D_{\alpha}$ is the generalized diffusion coefficient.

Processes known as moving boundary problems or Stefan like problems (in connection with the early work
of Joseph Stefan \cite{Ste91}) involve solution of the diffusion equation in a domain with a free 
boundary. They are governed by the formula (\ref{eq:diff}). Many monographs have been devoted 
to the classical Stefan problem, among which we can mention: \cite{Cra84,Hil87,Rub71,Gup03,Ozi93}.

The fractional Stefan problem is a generalization of the classical one, in which both the motion of the
molecules and the position of the moving boundary are governed by the low described by formula 
(\ref{eq:subdiff}). Currently, many analytical solutions (also approximate solutions) of the fractional 
Stefan problem are known and presented in papers 
\cite{Liu04,Vol10,Vol13,Sin11,Raj13,Xic08,Ros13,Ros16,Ala06,Liu09} and they constitute the 
overwhelming majority of all solutions. Numerical solutions and methods dedicated to solutions
are not so numerous \cite{Xia15,Vol14,Bla15,Bla14} and still require further research.

The numerical scheme discussed in this paper is an extension of the front fixing method developed in
\cite{Bla15}. The most important element of the new method is the much more accurate integration
algorithm for the fractional integro-differential equation. The new approach uses the Al-Alaoui 
operator, which significantly increases the accuracy of numerical integration. The optimal contribution
of the fractional trapezoidal rule and the fractional rectangular rule was implemented using an 
artificial neural network.

 \section{Preliminaries}
This section is devoted to integrals and derivatives of a non-integer order together with some 
of their properties. We also recall the two-parameter Wright function related to the theory of partial 
differential equations of fractional order. The definitions in the field of 
numerical methods will also be presented. At the beginning, let us focus our attention on the
definitions of two fractional operators: the left-sided Riemann-Liouville integral and the
left-sided Caputo derivative.
\begin{d1}
	The left-sided Riemann-Liouville integral of order $\alpha$, denoted as $I_{0+}^{\alpha}$, is given by the following formula for $Re(\alpha)>0$:
	\begin{equation}
	I_{0+}^{\alpha}f(t):=\frac{1}{\Gamma(\alpha)}\int_{0}^{t}\frac{f(u)du}
	{(t-u)^{1-\alpha}},
	\end{equation}
	where $\Gamma$ is the Euler gamma function.
\end{d1}

\begin{d2}
	Let $Re(\alpha)\in(0,1]$. The left-sided Caputo derivative of order $\alpha$
	is given by the formula:
	\begin{equation}
	{}^{c}D_{0+}^{\alpha}f(t):=\left\lbrace
	\begin{array}{ll}
	\frac{1}{\Gamma(1-\alpha)}
	\int_{0}^{t}\frac{f^{'}(u)du}{(t-u)^{\alpha}},&0<\alpha<1,\\
	\frac{d f(t)}{dt}, & \alpha=1.
	\end{array}
	\right.
	\label{eq:Cap}
	\end{equation}
\end{d2}
Two of the well known properties of integer-order integral and differential operators are
preserved by the following generalizations:
\begin{p1}[cf. Lemma~2.3 \cite{Kil06}]
	\label{pr:3}
	If $Re(\alpha)>0$, and $Re(\beta)>0$, then the equation
	\begin{equation}
	I_{0+}^{\alpha}I_{0+}^{\beta}f(t)=I_{0+}^{\alpha+\beta}f(t)
	\end{equation}
	is satisfied at almost every point $t\in[0,b]$ for $f(t)\in L_p(0,b)$ where 
	$1 \leq p\leq \infty$. If $\alpha+\beta>1$, then the above relation holds
	at any point of $[0,b]$.
\end{p1}

\begin{p1}[cf. Lemma~2.22, \cite{Kil06}]
	\label{pr:4}
	Let function $f\in C^1(0,b)$. Then, the composition rule for the left-sided Riemann-Liouville integral
	and the left-sided Caputo derivative is given as follows:
	\begin{equation}
	I_{0+}^{\alpha}{}^c D_{0+}^{\alpha}f(t)=f(t)-f(0).
	\end{equation}
\end{p1}
One of the very commonly used special functions in the theory of partial differential equations of
fractional order is two-parameter Wright function:
\begin{d2}
	Let $\gamma > -1$, $\delta\in \mathbb{C}, z \in \mathbb{C}$. The two-parameter Wright function
	 is given as the 
	following series:
	\begin{equation}
	\label{wright}
	W(z;\gamma,\delta):=\sum_{k=0}^{\infty}\frac{z^k}{k!\Gamma(\gamma k+\delta)},
	\end{equation}
\end{d2}
which is a generalization of the complementary error function.

The numerical scheme proposed in the further part of the paper uses a mesh of nodes defined
as follows:
\begin{d1}
	\label{def:grid}
	Let $\Gamma=\{(u,\tau): u\in (0,1); \tau\geq 0\}$ be a continuous region of solutions for
	the partial differential equation and $\tau^*=\frac{1}{p^{2/\alpha}}$ will be the end time. 
	Then the set 
	$\bar{\Gamma}=\{(u_i,\tau_j)\in \Gamma: x_i=i \Delta u, i\in\{0,1,...,m\}, \Delta u=
	\frac{1}{m}; \tau_j=j \Delta\tau, j\in\{ 0,1,...,n\};\Delta\tau=\frac{\tau^*}{n} \}$
	we call the rectangular regular mesh described by the set of nodes.
\end{d1}

The next three definitions will allow us to formulate the artificial neural network in the 
form of a very simple formula.
\begin{d1}
	\label{def:tanh}
	Function $\Tanh: R^n\to R^n$ such that\\
	$$\Tanh(x_1,x_2,...,x_n)= \left[ \tanh(x_1),\tanh(x_2),...,\tanh(x_n)\right], $$
	where $x_1,x_2,...,x_n\in R$ is a vector function, that assigns to each point 
	$[x_1,x_2,...,x_n]$ a n-dimensional vector lying in $R^n$ space.
\end{d1}
\begin{d1}
	\label{def:sig}
	Function $\Sig: R^n\to R^n$ such that\\
	$$\Sig(x_1,x_2,...,x_n)=
	 \left[ \frac{1}{1+e^{-x_1}},\frac{1}{1+e^{-x_2}},...,\frac{1}{1+e^{-x_n}}\right], $$
	where $x_1,x_2,...,x_n\in R$ is a vector function, that assigns to each point 
	$[x_1,x_2,...,x_n]$ a n-dimensional vector lying in $R^n$ space.
\end{d1}
\begin{d1}
	\label{def:psi}
	Function $\psi: R^n\to R^{n+1}$ such that\\
	$$\psi(x_1,x_2,...,x_n)=\left[x_1,x_2,...,x_n,1 \right], $$
	where $x_1,x_2,...,x_n\in R$ is a vector function, that assigns to each point 
	$[x_1,x_2,...,x_n]$ a n+1-dimensional vector lying in $R^{n+1}$ space.
\end{d1}
\section{Mathematical formulation of the problem}
\begin{figure}[h!t]
	\begin{center}
		\includegraphics[scale=1]{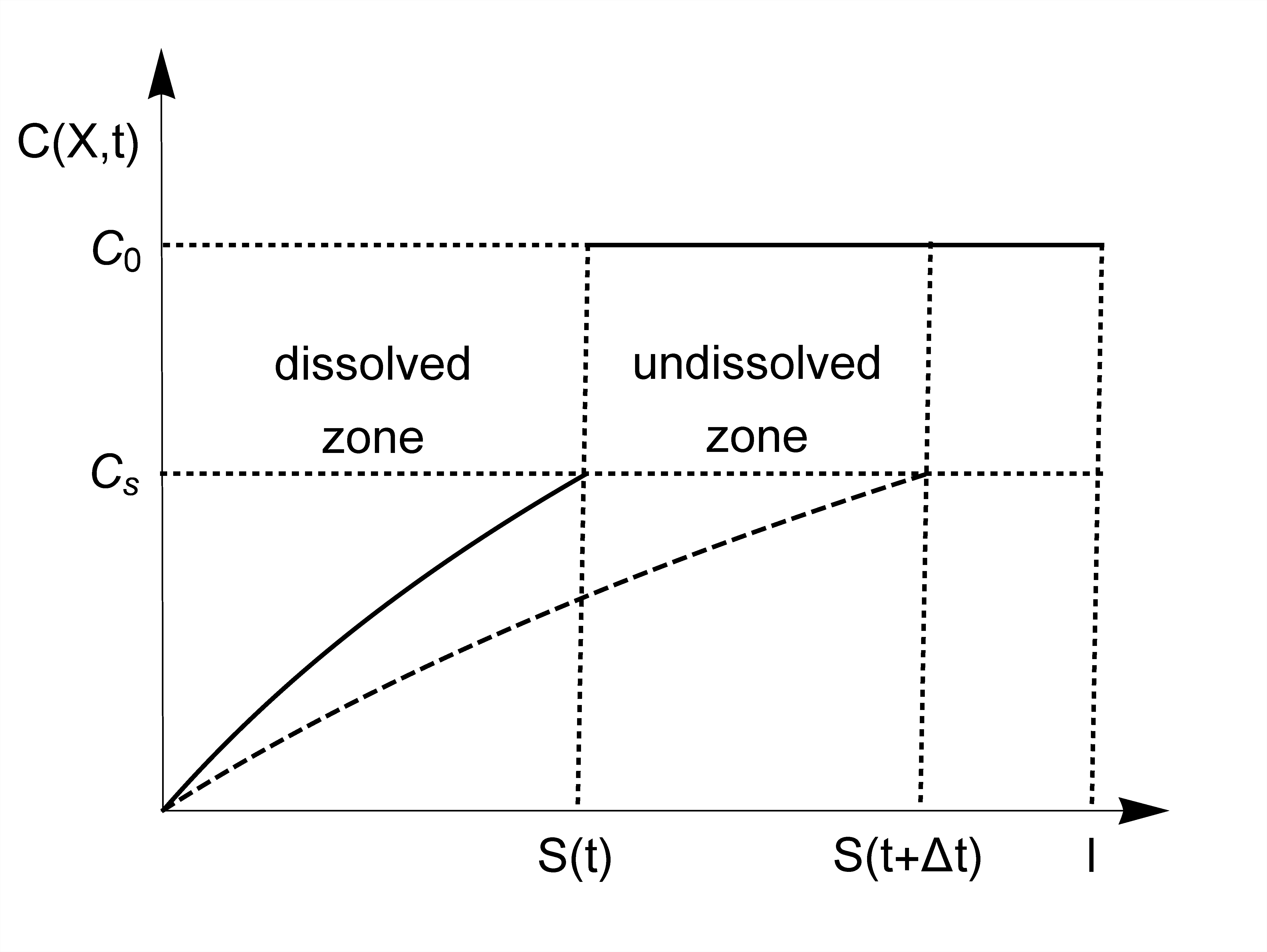}
		\caption{Slab dissolving from $x=0$ due to the zero drug concentration.}
		\label{fig:scheme}
	\end{center}
\end{figure}

The release of the solute is possible by anomalous diffusion through a penetrant fluid in the
polymer matrix consists of two regions: first bounded by $X=0$ and $X=S(t)$ where all solute is
dissolved and second one bounded by $X=S(t)$ and $X=l$ which contains undissolved solute. The
boundary of the domains formed by the diffusion front is described by the unknown function $S(t)$.
On the moving boundary $X = S(t)$ local mass conservation law is fulfilled and the concentration
$C$ is constant and equal to drug solubility $C_S$. We also assume that: $(i)$ generalized
diffusivity $\mathcal{D}_{\alpha}$ of the drug in the matrix is constant; $(ii)$ initial
concentration of drug $C_0$ is greater than its solubility $C_S$; $(iii)$ matrix is heterogeneous
and non-swellable. The mathematical model described diffusion of the solute through the dissolved
drug phase (illustrated by a simple scheme in Figure \ref{fig:scheme}) is formulated by
the subdiffusion equation:
\begin{equation}
\label{eq:sub}
{}^C D_{0+,t}^{\alpha}C(X,t)=\mathcal{D}_{\alpha}\frac{\partial^2 C(X,t)}{\partial X^2},
\quad 0<X<S(t), \quad t>0,
\end{equation}
supplemented with the boundary conditions
\begin{equation}
C(0,t)=0,\quad C(S(t),t)=C_S,\quad t>0,
\end{equation}
initial conditions
\begin{equation}
C(0,0)=C_S, \quad S(0)=0,
\end{equation}
and fractional Stefan condition
\begin{equation}
\label{eq:ste}
(C_0-C_S){}^C D_{0+,t}^{\alpha}C(X,t)=\mathcal{D}_{\alpha}\left.\frac{\partial C(X,t)}{\partial X}
\right|_{X=S(t)}.
\end{equation}
Model formulated by equations (\ref{eq:sub}-\ref{eq:ste}) assume that drug release is into a
perfect sink, with zero drug concentration i.e. $C(0, t) = 0$. This condition is a mathematical 
idealisation that is well approximated if the release medium is exchanged sufficiently rapidly to keep
sink conditions, or if the volume of the release medium is so large that drug concentration in the
medium is negligible.

Many mathematical models formulated by partial differential equations use dimensionless variables 
in numerical work. This approach is particularly useful in problems with universal scope 
applying to different scales \cite{Bar96,Bar03}. There are some advantages in using dimensionless
variables e.g. reduction of parameters in equations has a very positive effect on the time necessary
to carry out calculations for large sets of parameters.

Let us note that the above mathematical model depends on four parameters, where constant
generalized  diffusion coefficient $D_{\alpha}$ is measured in unit $cm^2\cdot s^{-\alpha}$.
Thus if we write the group $\tau=t\left( \frac{D_{\alpha}}{l^2}\right)^{1/\alpha}$, then we see that
$\tau$ is a dimensionless variable. For subdiffusion in a plane sheet we introduce the following
dimensionless variables:
\[
\tau=t\left( \frac{D_{\alpha}}{l^2}\right)^{1/\alpha},\quad x=\frac{X}{l}, \quad 
s(\tau)=\frac{S(t)}{l}, \quad c=\frac{C}{C_S},
\]
which reduces the number of parameters to two.
In our case, a whole set of solutions for different values of: generalized  diffusion coefficient 
$D_{\alpha}$,  initial concentration of drug $C_0$, drug solubility $C_S$ and matrix thickness $l$
can be obtained from a single solution in dimensionless variables by simple scaling. So we 
can rewrite  equations (\ref{eq:sub}-\ref{eq:ste}) in a dimensionless form:

\begin{equation}
\label{eq:nonsub}
{}^C D_{0+,\tau}^{\alpha}c(x,\tau)=\frac{\partial^2 c(x,\tau)}{\partial x^2},
\quad 0<x<s(\tau), \quad \tau>0,
\end{equation}
supplemented with the boundary conditions
\begin{equation}
c(0,\tau)=0,\quad c(s(\tau),\tau)=1,\quad \tau>0,
\end{equation}
initial conditions
\begin{equation}
c(0,0)=1, \quad s(0)=0,
\end{equation}
and fractional Stefan condition
\begin{equation}
\label{eq:nonste}
{}^C D_{0+,\tau}^{\alpha}c(x,\tau)=\lambda \left.\frac{\partial c(x,\tau)}{\partial x}
\right|_{x=s(\tau)},
\end{equation}
where $\lambda=\left(\frac{C_0}{C_S} \right)^{-1} $ is fractional Stefan number.

From the point of view of further numerical results presented in this paper, the analytical solution
of the one-phase fractional Stefan problem obtained by Liu Junyi et al. \cite{Liu09} is very valuable,
because it will allow validation of the method proposed in section four. Using the notation 
introduced in this paper, the closed analytical solution to the problem under consideration is
expressed in terms of the Wright’s function:
\begin{equation}
\label{eq:exa}
c(x,\tau)=\frac{1-W\left(-\frac{x}{\tau^{\frac{\alpha}{2}}};-\frac{\alpha}{2},
	1 \right) }{1-W\left(-p;-\frac{\alpha}{2},1 \right)},
\end{equation}
where $p$ is the coefficient of the power function $s(\tau)=p\tau^{\frac{\alpha}{2}}$ 
describing the position of the diffusion front which can be determined from the transcendental 
equation
\begin{equation}
\label{eq:trans}
\frac{\lambda \Gamma\left(1-\frac{\alpha}{2} \right) }{p \Gamma\left(1+
	\frac{\alpha}{2} \right) }=\frac{1-W\left(-p;-\frac{\alpha}{2},1 \right)}
	{W\left(-p;-\frac{\alpha}{2},1-\frac{\alpha}{2} \right)}.
\end{equation}
\section{Solution of the problem}
The subject of our considerations is to find a numerical solution of the problem described by 
the equations (\ref{eq:sub}-\ref{eq:ste}). The proposed approach takes place in three stages:
\begin{enumerate}
	\item we introduce to subdiffusion equation (\ref{eq:nonsub}) new spatial variable,
	which immobilizes dissolution boundary for all $\tau$. Next we convert the 
	fractional partial differential equation to the equivalent fractional integro-differential
	equation. We are using two methods of numerical integration: the fractional rectangle rule and
	the fractional trapezoidal rule to discretization of the integro-differential equation. The 
	share of the trapezoidal rule in the integration scheme is determined by the 
	parameter $\phi\in[0,1]$, whose optimal value at this stage is unknown.
	\item in the second stage, we determine the optimal value of the $\phi$ parameter. This 
	task is implemented as follows: for fixed values of the model and grid parameters (the value
	of $p$ we get from the equation (\ref{eq:trans})) we minimize some functional expressing
	the sum of the absolute errors generated by the numerical scheme received in the previous
	point. The $\phi$ values obtained in this way and the associated values of the model and
	mesh parameters create training set for the artificial neural network. Next, we create the
	artificial neural network to predict the values of the $\phi$ parameter.
	\item At this stage, the value of the $\phi$ parameter is already known and can be used in
	the scheme from the first point. To make a numerical scheme independent of the analytical 
	solution, we formulate an iterative algorithm that allows us to determine the value of the 
	parameter $p$ regardless of the equation (\ref{eq:trans}).
	
\end{enumerate}
\subsection{Front fixing method}
The basic idea of the front fixing method is a suitable choice of a new space variable for
the governing equation. The transformation:
\begin{equation}
\label{eq:sim}
u=\frac{x}{p \tau^\frac{\alpha}{2}},
\end{equation}
fixes dissolution boundary at the point $u=1$ for all $\tau$. By using the following
 relationships 
\begin{equation}
\frac{\partial^2 c(x,\tau)}{\partial x^2}=\frac{1}{p^2 \tau^{\alpha}}\frac{\partial^2 c(u,\tau)}{\partial u^2}
\end{equation}
\begin{equation}
\label{eq:caputo}
\begin{split}
{}^C D_{0+,\tau}^{\alpha}c(u,\tau)=&
\frac{1}{\Gamma(1-\alpha)}\int_{0}^{\tau}\frac{1}{(\tau-\xi)^{\alpha}}
\frac{\partial c(u,\xi)}{\partial \xi} d\xi-\\&
\frac{\alpha u}{2\Gamma(1-\alpha)}\int_{0}^{\tau}\frac{1}{(\tau-\xi)^{\alpha}}\frac{1}{\xi}
\frac{\partial c(u,\xi)}{\partial u} d\xi
\end{split}
\end{equation}
and Property \ref{pr:3} and \ref{pr:4}
\begin{equation}
\label{eq:sub3}
\begin{split}
I_{0+,\tau}^{\alpha}{}^{C}D_{0+,\tau}^{\alpha}c(u,\tau)-\frac{\alpha 
	u}{2}I_{0+,\tau}^{\alpha} I_{0+,\tau}^{1-\alpha}\left(\frac{\partial}{\partial u}
\frac{c(u,\tau)}{\tau}\right)=\frac{1}{p^2}I_{0+,\tau}^{\alpha}
\left(\frac{1}{\tau^{\alpha}}\frac{\partial^2 c(u,\tau)}{\partial u^2}\right),
\end{split}
\end{equation}
the subdiffusion equation (\ref{eq:nonsub}) can be written
\begin{equation}
\label{eq:fin1}
\begin{split}
c(u,\tau)&=c(u,0)+\frac{\alpha u}{2}\int_{0}^{\tau}
\frac{\partial}{\partial u}\frac{c(u,\xi)}{\xi}d\xi+
\frac{1}{p^2 \Gamma(\alpha)}\int_{0}^{\tau}\frac{1}{(\tau-\xi)^{1-\alpha}}
\frac{1}{\xi^{\alpha}}\frac{\partial^2 c(u,\xi)}{\partial u^2}d\xi.
\end{split}
\end{equation}
Now, we introduce the auxiliary function
\begin{equation}
\bar{c}(u,\tau)=c(u,\tau)\tau^{-\alpha}.
\end{equation}
So, finally, we get the integro-differential equation
\begin{equation}
\label{eq:fin2}
\begin{split}
&\bar{c}(u,\tau)\tau^{\alpha}=\bar{c}(u,0)\tau_{0}^{\alpha}+\frac{\alpha u}{2}\int_{0}^{\tau}
\frac{\partial}{\partial u}\frac{\bar{c}(u,\xi)}{\xi^{1-\alpha}}d\xi+\\&+
\frac{1}{p^2 \Gamma(\alpha)}\int_{0}^{\tau}\frac{1}{(\tau-\xi)^{1-\alpha}}
\frac{\partial^2 \bar{c}(u,\xi)}{\partial u^2}d\xi,
\end{split}
\end{equation}
supplemented by the boundary conditions
\begin{equation}
\bar{c}(0,\tau)=0,\quad \bar{c}(1,\tau)=\tau^{-\alpha},
\end{equation}
and initial condition
\begin{equation}
\bar{c}(u,0)=\lim_{\tau_0\to 0^{+}} \tau_{0}^{-\alpha}.
\end{equation}
It should be noted that due to further numerical considerations, we assume that $\tau_0$ is
a very small positive real number.

Here we should explain our new approach. We can see that the integro-differential equation
discussed above, and in particular the last integral term on the right hand side of this
equation can be discretized in different ways. In \cite{Bla15}, the fractional trapezoidal rule 
was used to approximate the integral. However, it turns out that the use of the Al-Alaoui 
operator (\cite{Ala93},\cite{Ala06},\cite{Ala08}) significantly improves the integration 
accuracy, which will be shown in examples later in the paper.

The second component on the right hand side of equation (\ref{eq:fin2}) can be approximated using
the central differential quotient to discretize the derivative of function $c$ with respect 
to variable $u$ and the rectangle rule for calculating the integral:
\begin{equation}
\label{eq:in1}
\begin{split}
&\frac{\alpha u_i}{2}\int_{0}^{\tau_{k+1}} \left( \frac{\partial \bar{c}(u,\xi)}
{\partial u}\right)_{i,j} \frac{1}
{\xi^{1-\alpha}}d\xi\approx \sum_{j=1}^{k+1}q_{i,j}(\bar{c}_{i+1,j}-\bar{c}_{i-1,j}),
\end{split}
\end{equation}
where
\begin{equation}
q_{i,j}=\frac{\alpha i \tau_j^{\alpha-1}\Delta \tau}{4}.
\end{equation}
Now let us apply the integration concept proposed by Al-Alaoui. First, we are using the
fractional trapezoidal rule for the left-sided Riemann-Liouville integral and the differential 
quotient to discretize the second derivative
\begin{equation}
\begin{split}
&\frac{1}{p^2 \Gamma(\alpha)}\int_{0}^{\tau_{k+1}}\frac{1}{(\tau_{k+1}-\xi)^{1-\alpha}}
\left( \frac{\partial^2 \bar{c}(u,\xi)}{\partial u^2}\right)_{i,j} d\xi\approx 
\sum_{j=0}^{k+1}r_{j,k+1}(\bar{c}_{i+1,j}-2\bar{c}_{i,j}+\bar{c}_{i-1,j}),
\end{split}
\end{equation}
where
\begin{equation}
\label{eq:weight}
r_{j,k+1}=\frac{1}{p^2 \Gamma(\alpha)(\Delta u)^2}\left\{
\begin{array}{ll}
\frac{(\Delta\tau)^{\alpha}}{\alpha(\alpha+1)}\left( k^{\alpha+1}-(k-\alpha)(k+1)^{\alpha}
\right) & \textrm{for $j = 0$} \\
\frac{(\Delta\tau)^{\alpha}}{\alpha(\alpha+1)}( (k-j+2)^{\alpha+1}+(k-j)^{\alpha+1}&\\
-2(k-j+1)^{\alpha+1}) & \textrm{for $ 1\leq j \leq k$}\\
\frac{(\Delta\tau)^{\alpha}}{\alpha(\alpha+1)}& \textrm{for $j=k+1$}
\end{array} \right.
\end{equation}
Then we calculate the same integral, but this time using the fractional rectangle rule
\begin{equation}
\begin{split}
&\frac{1}{p^2 \Gamma(\alpha)}\int_{0}^{\tau_{k+1}}\frac{1}{(\tau_{k+1}-\xi)^{1-\alpha}}
\left( \frac{\partial^2 \bar{c}(u,\xi)}{\partial u^2}\right)_{i,j} d\xi\approx 
\sum_{j=0}^{k}w_{j,k+1}(\bar{c}_{i+1,j}-2\bar{c}_{i,j}+\bar{c}_{i-1,j}),
\end{split}
\end{equation}
where
\begin{equation}
w_{j,k+1}=\frac{1}{p^2 \Gamma(\alpha)} \int_{\tau_j}^{\tau_{j+1}}\frac{1}
{(\tau_{k+1}-\xi)^{1-\alpha}}d\xi=\frac{(\Delta \tau)^{\alpha}}{p^2 
\Gamma(\alpha+1)}((k+1-j)^{\alpha}-(k-j)^{\alpha}).
\end{equation}
Finally, the approximation of the left-sided Riemann-Liouville integral can be expressed 
in the following form
\begin{equation}
\label{eq:in2}
\begin{split}
&\frac{1}{p^2 \Gamma(\alpha)}\int_{0}^{\tau_{k+1}}\frac{1}{(\tau_{k+1}-\xi)^{1-\alpha}}
\left( \frac{\partial^2 \bar{c}(u,\xi)}{\partial u^2}\right)_{i,j} d\xi\approx\\ \approx&
\phi r_{k+1,k+1}(\bar{c}_{i+1,k+1}-2\bar{c}_{i,k+1}+\bar{c}_{i-1,k+1})+\\+&
\sum_{j=0}^{k}(\phi r_{j,k+1}+(1-\phi) w_{j,k+1})(\bar{c}_{i+1,j}-2\bar{c}_{i,j}+\bar{c}_{i-1,j})
\end{split}
\end{equation}
Using formulas (\ref{eq:in1}-\ref{eq:in2}) and (\ref{eq:fin2}), we get a numerical scheme that
can be written in the matrix form
\begin{equation}
\label{eq:schemat}
\mathbf{A}\mathbf{\bar{C}}_{k+1}=\mathbf{B},
\end{equation}
where \textbf{A} is an $(m-1)\times (n-1)$ matrix of coefficients of system of linear equations
\[\scalebox{0.9}{$\mathbf{A}=\left[\begin{array}{cccccccc}
a_{k+1}^2 & a_{1,k+1}^3 & 0 & 0 & \cdots & 0 & 0 & 0 \\
a_{2,k+1}^1 & a_{k+1}^2 & a_{2,k+1}^3 & 0 & \cdots & 0 & 0 & 0 \\
0 & a_{3,k+1}^1 & a_{k+1}^2 & a_{3,k+1}^3 & \cdots & 0 & 0 & 0 \\
\vdots & \vdots & \vdots & \vdots & \ddots & \vdots & \vdots & \vdots\\
0 & 0 & 0 & a_{i,k+1}^1 & a_{k+1}^2 & a_{i,k+1}^3 & 0 & 0\\
\vdots & \vdots & \vdots & \vdots & \ddots & \vdots & \vdots & \vdots\\
0 & 0 & 0 & 0 & \cdots & a_{m-2,k+1}^1 & a_{k+1}^2 & a_{m-2,k+1}^3 \\
0 & 0 & 0 & 0 & \cdots & 0 & a_{m-1,k+1}^1 & a_{k+1}^2 \\
\end{array}
\right]_{(m-1)\times(m-1)},
$}\]
\textbf{B} is a column vector with $n-1$ entires
$$\mathbf{B}=\left[\begin{array}{c}
b_1-a_{1,k+1}^1 \bar{c}_{0,k+1}\\
b_2\\
b_3\\
\vdots\\
b_i\\
\vdots\\
b_{m-2}\\
b_{m-1}-a_{m-1,k+1}^3 \bar{c}_{m,k+1}
\end{array}
\right]_{(m-1)\times 1},
$$
and $\mathbf{\bar{C}}_{k+1}$ is a column vector with $m-1$ entires. We define the elements
of matrix \textbf{A} and \textbf{B} as follows:
\begin{eqnarray*}
	&&a_{i,k+1}^1:=-\phi r_{k+1,k+1}+q_{i,k+1},\\
	&&a_{k+1}^2:=\tau_{k+1}^{\alpha}+2\phi r_{k+1,k+1},\\
	&&a_{i,k+1}^3:=-\phi r_{k+1,k+1}-q_{i,k+1},\\
	&&b_i:=\bar{c}_{i,0}\tau_0^{\alpha}+\sum_{j=0}^{k}
	(\phi r_{j,k+1}+(1-\phi) w_{j,k+1})(\bar{c}_{i+1,j}-2\bar{c}_{i,j}+\bar{c}_{i-1,j})+\\
	&&+\sum_{j=1}^{k}q_{i,j} (\bar{c}_{i+1,j}-\bar{c}_{i-1,j}).
\end{eqnarray*}
The formulas below allow us to recover the original spatial variable $x$ and function $c$.
\begin{equation}
c_{i,j}=\bar{c}_{i,j}\tau_{j}^{\alpha},
\end{equation}
\begin{equation}
x_{i,j}=u_{i} p \tau_{j}^{\alpha/2}.
\end{equation}
Particular attention should be given to the fact that above-derived numerical scheme depends
on two unknown parameters: $p$ and $\phi$. In the next two subsections we will demonstrate how to 
determine their values.
\subsection{Prediction of the $\phi$ parameter value using an artificial neural network}
Simultaneous determination of the optimum value of the parameter $\phi$ and $p$ is very
problematic. Therefore, this task will be divided into two stages. First, we determine the
optimal $\phi$ value. We will use an artificial neural network for this purpose, because this
 mathematical construct can generate predictions for very complex problems. Let us assume
 that the parameter $\phi$ depends on three variables: $\alpha$, $\lambda$ and $\Delta u$ 
 (second grid parameter $\Delta\tau$ depends on the spatial step as follows- in definition
 \ref{def:grid} adopted $n=4m$). The proposed network is a four-layer, unidirectional and a
 full neural network. The first network layer (input layer) consists of three neurons loading
 variables: $\Delta u, \lambda, \alpha$ and bias neuron with a fixed value of 1. The 
 first hidden layer of the network contains five neurons applying $\tanh$ as an activation 
function and bias neuron with a fixed value of 1. The next hidden layer also contains five 
neurons using a sigmoid function as an activation function and bias neuron with a fixed value 
of 1. The output layer contains one neuron applying sigmoid function as an activation function.
Two criteria were used to select the number and type of neurons: the accuracy of the results
generated by the network and it's as simple structure as possible. As we know, a large number
of neurons in the hidden layer increases the computing power and flexibility of the neural
network while learning complicated patterns, but also contributes to good matching of the 
network to the training set and loss of generalizing ability on the test set. The learning 
process (adjustment of the weights to every connection between neurons) of the neural network 
was carried out using a backpropagation algorithm.

The neural network learning process requires a training set, which is a group of sample 
inputs we can feed into the neural network in order to train the model. The training set used
in this paper was created as follows. We first determine the possible values of loading
variables: $\Delta u, \lambda, \alpha$:
\begin{eqnarray*}
	&& \Delta u\in z_1=\left\lbrace \frac{1}{25},\frac{1}{50},\frac{1}{75},\frac{1}{100}
	\right\rbrace, \\
	&& \lambda\in z_2=\left\lbrace0.25, 0.5, 0.75, 1, 1.25, 1.5, 1.75, 2,
	2.25, 2.5, 2.75, 3 \right\rbrace,\\
	&& \alpha\in z_3=\left\lbrace 0.2, 0.4, 0.6, 0.8, 1\right\rbrace. 
\end{eqnarray*}
 Next, we solve a so-called inverse problem for set $\{z_1\times z_2 \times z_3\}$. In an 
 inverse problem, some information is not known, in our case we do not know parameter $\phi$.
 By minimizing the following functional
\begin{equation}
F(\phi)=\sum_{j=0}^{n}\sum_{i=0}^{m} \left| c_{i,j}(\phi)-c(x_i,\tau_j)\right| ,
\end{equation}
we find the approximate value of the parameter $\phi.$ By $c_{i,j}(\phi)$ and $c(x_i,\tau_j)$
we denoted numerical solution obtained in mesh node $(i \Delta x, j \Delta \tau)$ by the scheme 
proposed in subsection 4.1, where the value of $p$ is known- obtained from formula 
(\ref{eq:trans}) and closed analytical solution obtained from formulas (\ref{eq:exa}),
 (\ref{eq:trans}) respectively.

\begin{figure}[h!t]
	\begin{center}
		\includegraphics[scale=0.8]{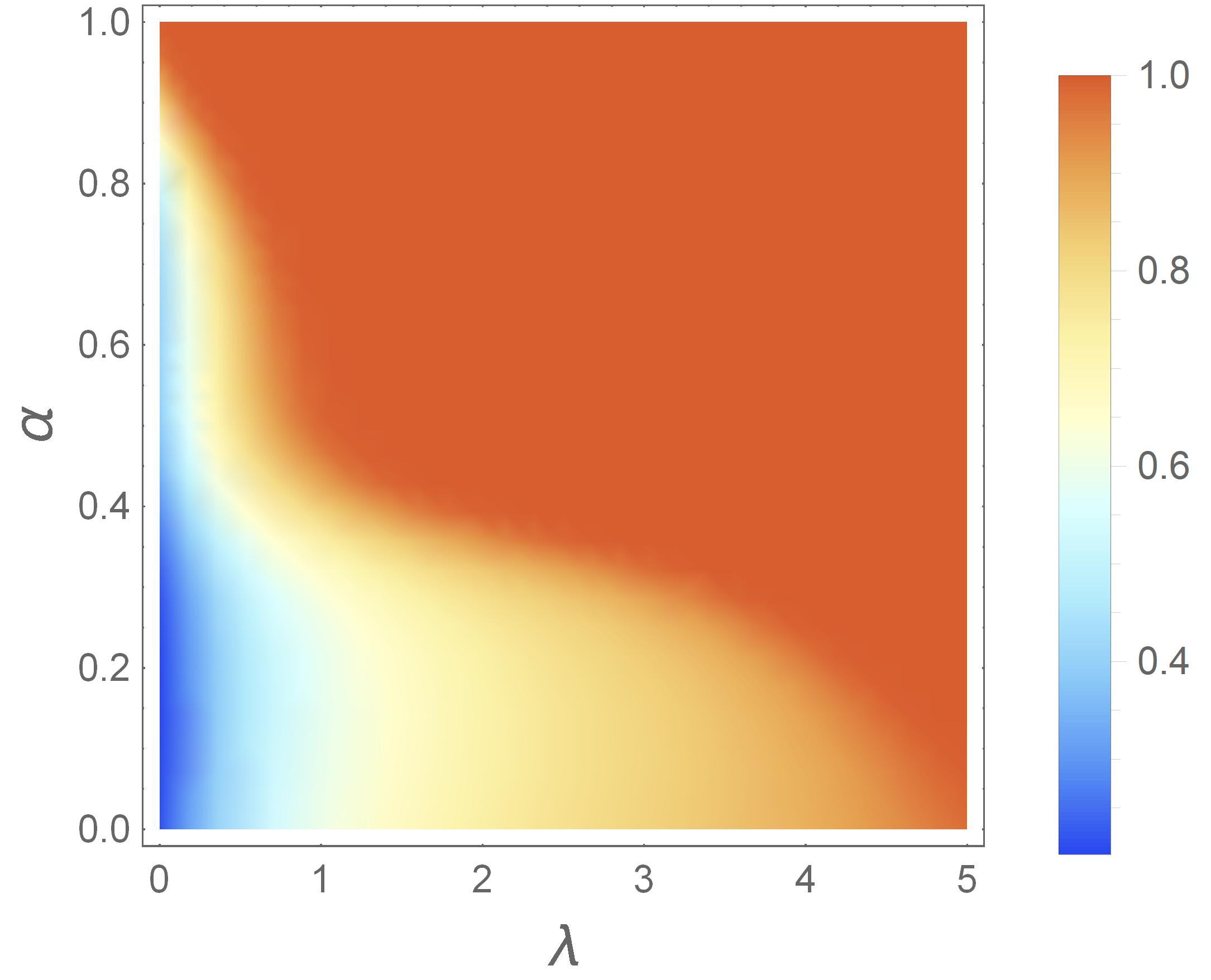}
		\caption{Plot of the $\phi$ function for $\alpha\in(0,1)$, $\lambda\in (0,5)$
		and $\Delta x=0.005$.}
		\label{fig:0}
	\end{center}
\end{figure}

As a result of the neural network learning process, a set of weights was obtained that we 
wrote in the matrix form:
$$
W_1:=\left(
\begin{array}{ccccc}
	6.02145 & 2.14126 & 0.677927 & -13.8801 & 1.52535 \\
	-0.514316 & 0.150222 & -0.853299 & 0.162208 & -0.88522 \\
	0.80747 & -7.26578 & -0.0383236 & -0.307713 & -5.98303 \\
	-0.140854 & 2.3685 & 0.0361885 & -0.284572 & 5.68518 \\
\end{array}
\right),
$$
$$
W_2:=\left(
\begin{array}{ccccc}
-7.18748 & -2.4983 & -0.72953 & -1.40697 & 0.748731 \\
-8.66802 & -2.94062 & 0.0381687 & -0.909638 & 0.268698 \\
-1.77908 & 0.0049326 & -4.28149 & -0.264405 & 0.382544 \\
11.2767 & 2.75494 & -1.16368 & -2.87708 & -1.78879 \\
-7.43406 & -1.53303 & -1.25246 & -1.20859 & 1.33437 \\
-0.590041 & -0.72953 & 1.50501 & -2.31013 & -0.662821 \\
\end{array}
\right),
$$
$$
W_3:=\left(
\begin{array}{c}
10.4966 \\
3.34645 \\
3.03464 \\
4.17581 \\
-3.59841 \\
-0.490178 \\
\end{array}
\right).
$$
The $W_1$ matrix contains the weight values connecting the network input layer with the 
first hidden layer. The row number of matrix $W_1$ means the number of the neuron in the 
input layer, while the number of the matrix column means the number of the neuron in the 
first hidden layer. Matrixes $M_2$ and $M_3$ have an analogous interpretation as in the case
of the $M_1$ matrix, but refer to a subsequent network layers. Using the above matrix notation,
we can write our neural network in the form of the following function of three variables:
\begin{equation}
\phi(\Delta u,\lambda, \alpha):=\Sig(\psi(\Sig(\psi(\Tanh([\Delta u,\lambda, \alpha,1]
\times W_1))\times W_2))\times W_3),
\end{equation}
whose values are illustrated in the Figure \ref{fig:0}.

An important step in creating artificial neural networks is the process of validation, in 
which we use another group of sample inputs which were not included in training and are 
different from the samples in the training set. The 5th section of this paper is devoted 
to the validation process of the proposed artificial neural network.
\subsection{Determination of the parameter $p$}
A boundary surface, on which dissolution occurs, moves across the slab separating a region of
a dissolved drug from an undissolved core as in Figure \ref{fig:scheme}. The position of the 
dissolution front is described by the function $s$ (dimensionless form) whose two values 
are known: at the beginning of the dissolution process $s=0$ and when dissolution is complete
$s=1$. Based on this simple observation, we formulate the convergence criterion
\begin{equation}
\label{eq:conver}
|1-p\tau_n^{\frac{\alpha}{2}}|<\epsilon,
\end{equation}
where
\begin{equation}
\label{eq:p}
p=\frac{1}{n+1}\sum_{j=0}^{n}\sqrt{\lambda \frac{\Gamma(1-\frac{\alpha}{2})\Delta c_j}
	{\Gamma(1+\frac{\alpha}{2})\Delta u}}.
\end{equation}
The proposed algorithm is a simple and more elegant version of the algorithm from \cite{Bla15},
that can be written on the following block diagram
\begin{figure}[h!t]
	\begin{center}
		\includegraphics[scale=0.45]{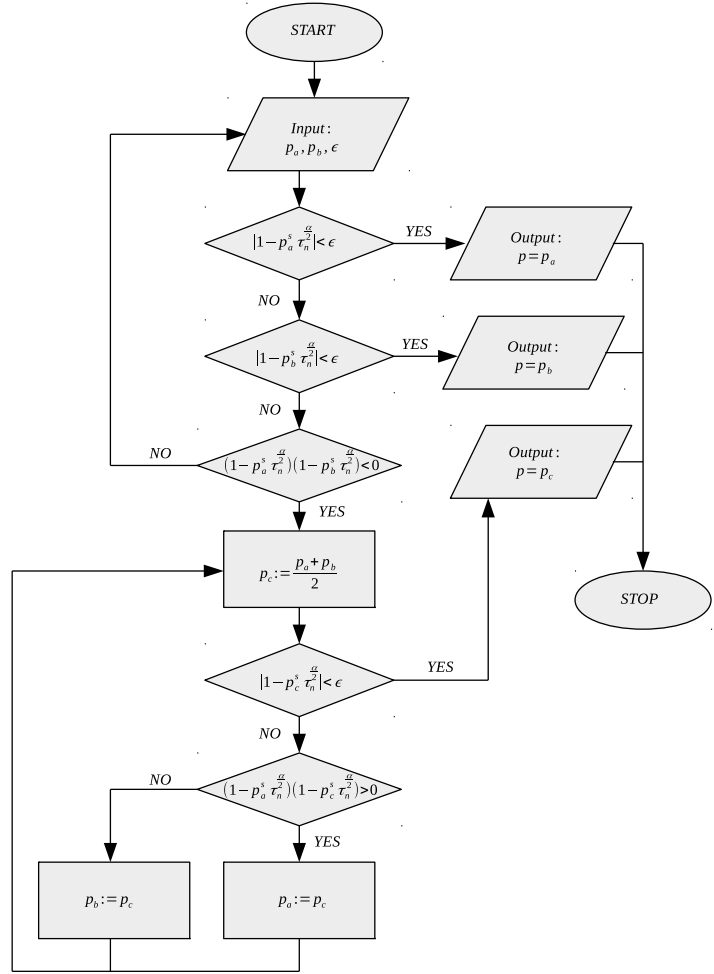}
		\caption{A block diagram of an iterative algorithm for determining the 
			coefficient $p$}
		\label{fig:sb}
	\end{center}
\end{figure}

The block diagram presented in Figure \ref{fig:sb} requires some additional detailed
explanations, special attention should be given to the superscript $s$ of the coefficient $p$.
The difference in notation has its interpretation. The values: $p_a$, $p_b$, $p_c$ contribute to
the formation of a new spatial variable (\ref{eq:sim}) used in numerical scheme
(\ref{eq:schemat})-they allow us to determine the end times $\tau_n$, while the superscript $s$
means that the values: $p_a^s$, $p_b^s$, $p_c^s$ have been calculated using the formula 
(\ref{eq:p}). The iterative algorithm ends when $p_{*}$ and $p_*^s$ are almost identical,
when condition (\ref{eq:conver}) is fulfilled.\newpage

\section{Numerical examples}
\begin{figure}[h!t]
	\begin{center}
		\includegraphics[scale=0.88]{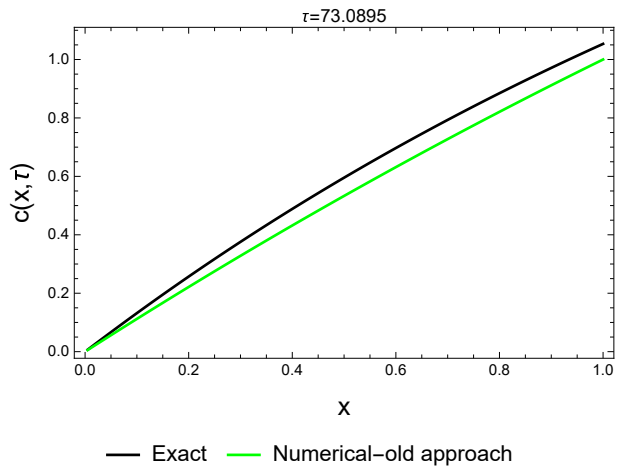}
		\caption{Numerical solution of function $c$ for $\lambda=\frac{1}{3}$ and $\alpha=0.25$.}
		\label{fig:1}
	\end{center}
\end{figure}

\begin{figure}[h!t]
	\begin{center}
		\includegraphics[scale=0.88]{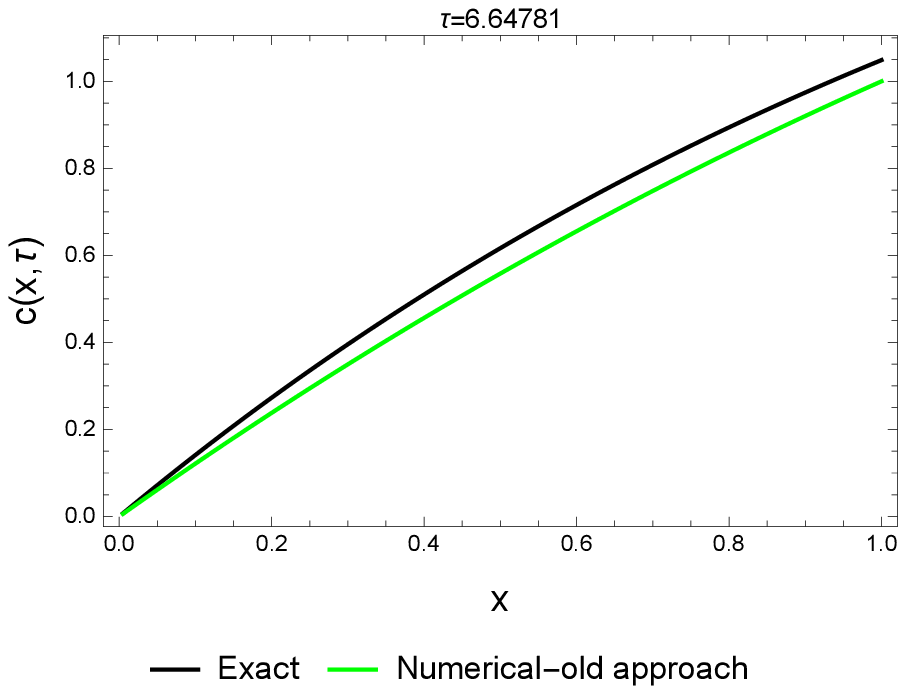}
		\caption{Numerical solution of function $c$ for $\lambda=\frac{2}{3}$ and $\alpha=0.25$.}
		\label{fig:6}
	\end{center}
\end{figure}

\begin{figure}[h!t]
	\begin{center}
		\includegraphics[scale=0.88]{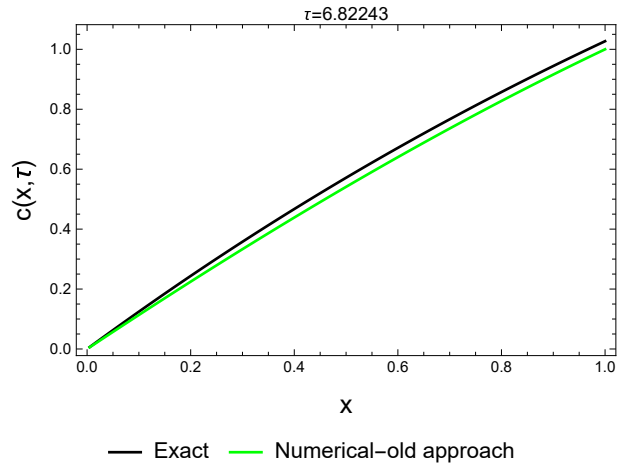}
		\caption{Numerical solution of function $c$ for $\lambda=\frac{1}{3}$ and $\alpha=0.5$.}
		\label{fig:11}
	\end{center}
\end{figure}

\begin{figure}[h!t]
	\begin{center}
		\includegraphics[scale=0.88]{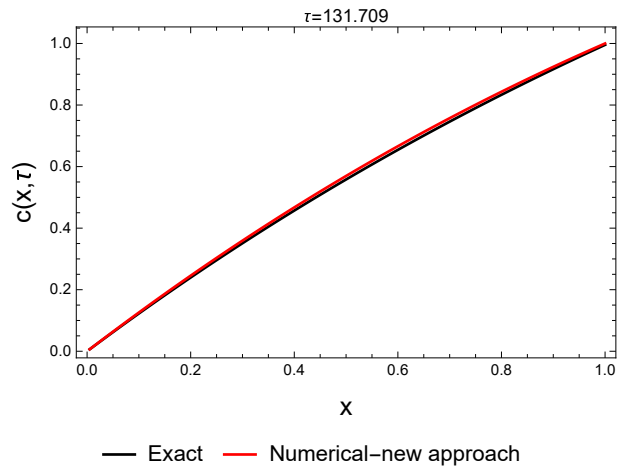}
		\caption{Numerical solution of function $c$ for $\lambda=\frac{1}{3}$ and $\alpha=0.25$.}
		\label{fig:2}
	\end{center}
\end{figure}

\begin{figure}[h!t]
	\begin{center}
		\includegraphics[scale=0.88]{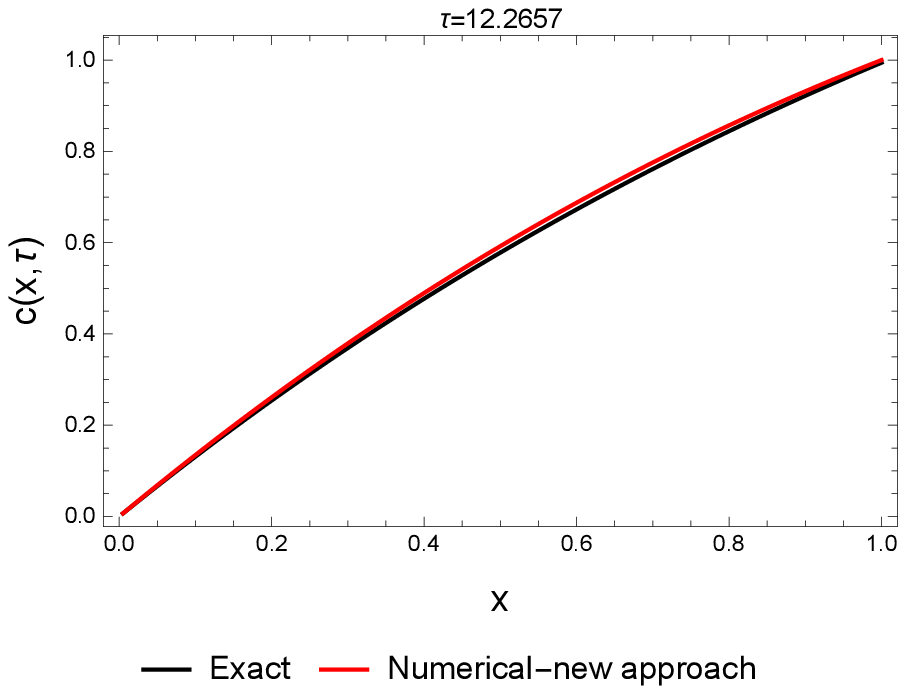}
		\caption{Numerical solution of function $c$ for $\lambda=\frac{2}{3}$ and $\alpha=0.25$.}
		\label{fig:7}
	\end{center}
\end{figure}

\begin{figure}[h!t]
	\begin{center}
		\includegraphics[scale=0.88]{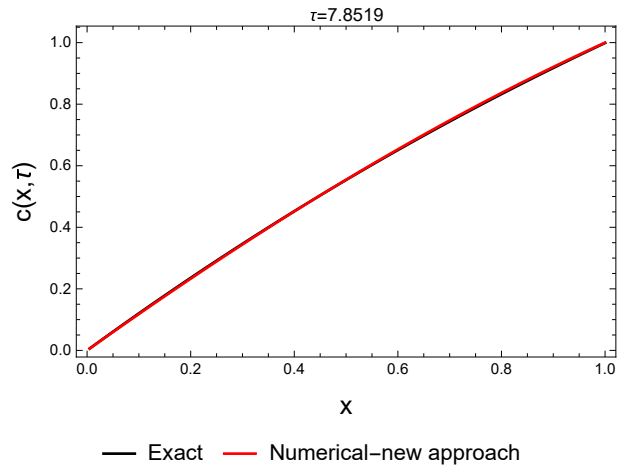}
		\caption{Numerical solution of function $c$ for $\lambda=\frac{1}{3}$ and $\alpha=0.5$.}
		\label{fig:12}
	\end{center}
\end{figure}

\begin{figure}[h!t]
	\begin{center}
		\includegraphics[scale=0.88]{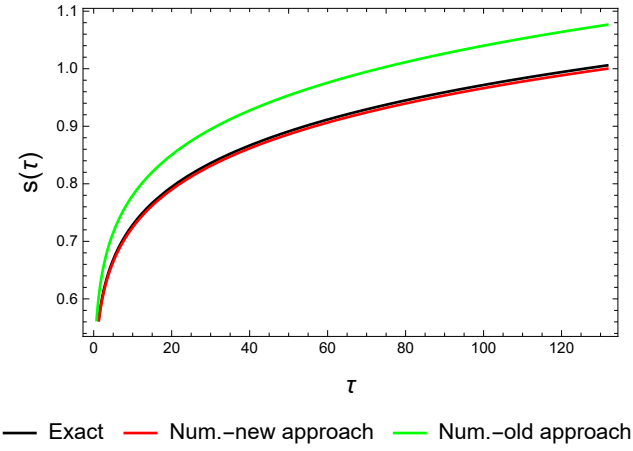}
		\caption{Numerical solution of function $s$ for $\lambda=\frac{1}{3}$ and $\alpha=0.25$.}
		\label{fig:3}
	\end{center}
\end{figure}

\begin{figure}[h!t]
	\begin{center}
		\includegraphics[scale=0.88]{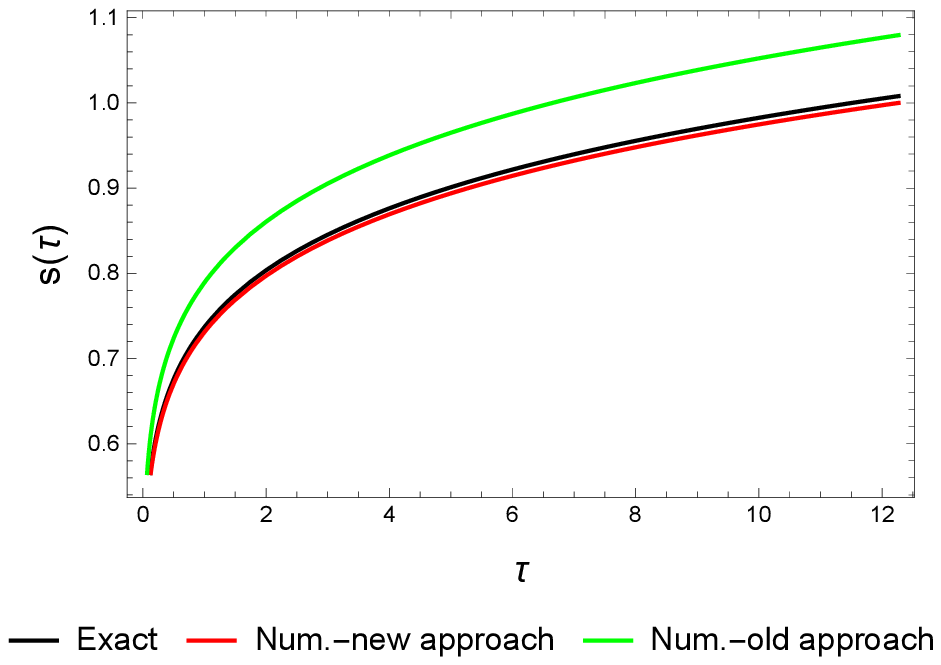}
		\caption{Numerical solution of function $s$ for $\lambda=\frac{2}{3}$ and $\alpha=0.25$.}
		\label{fig:8}
	\end{center}
\end{figure}

\begin{figure}[h!t]
	\begin{center}
		\includegraphics[scale=0.88]{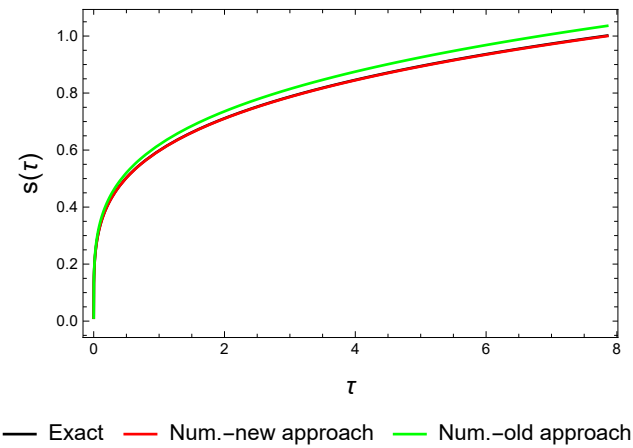}
		\caption{Numerical solution of function $s$ for $\lambda=\frac{1}{3}$ and $\alpha=0.5$.}
		\label{fig:13}
	\end{center}
\end{figure}

\begin{figure}[h!t]
	\begin{center}
		\includegraphics[scale=0.75]{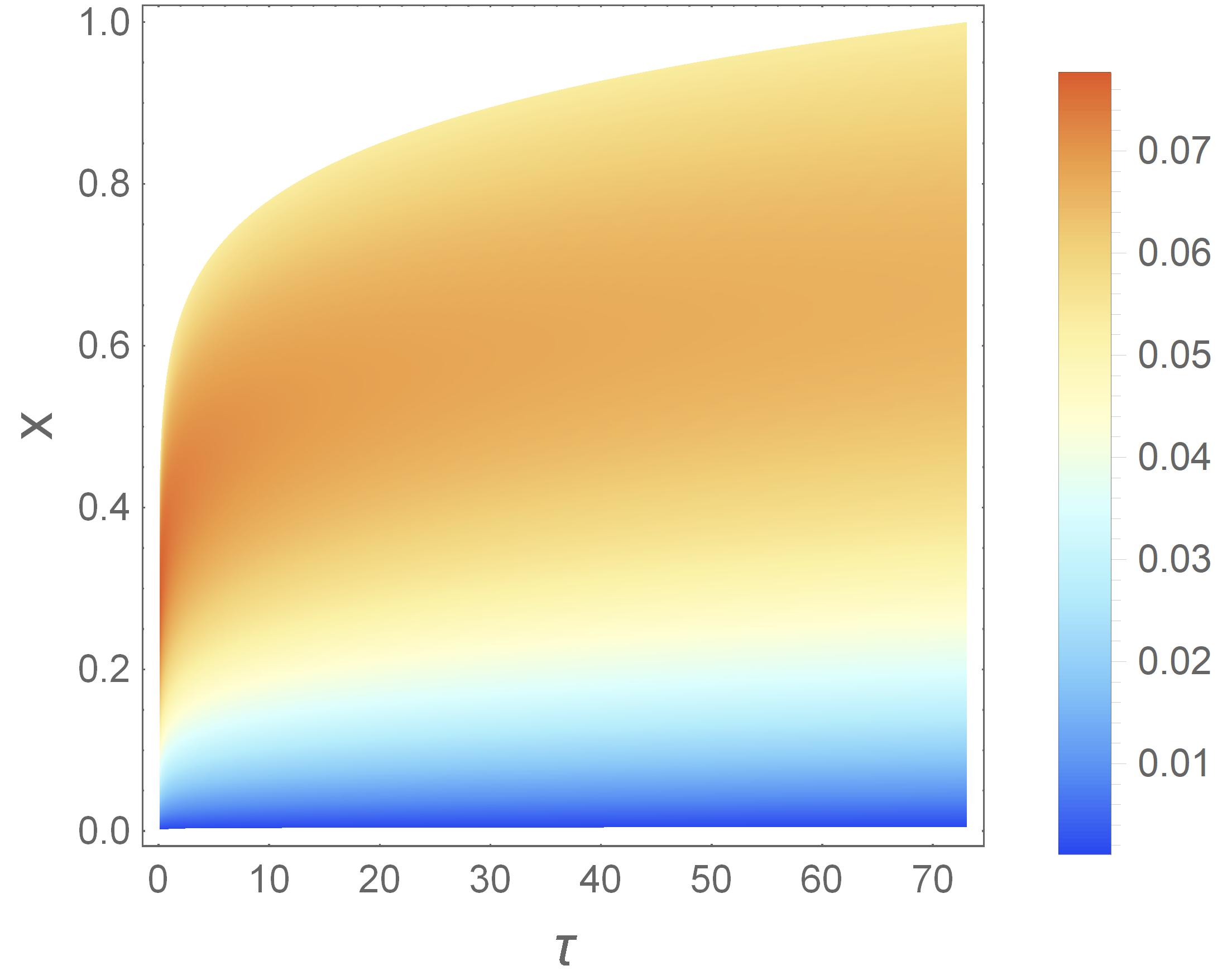}
		\caption{Absolute error generated by the old version of the numerical scheme for 
		$\lambda=\frac{1}{3}$ and $\alpha=0.25$.}
		\label{fig:4}
	\end{center}
\end{figure}

\begin{figure}[h!t]
	\begin{center}
		\includegraphics[scale=0.75]{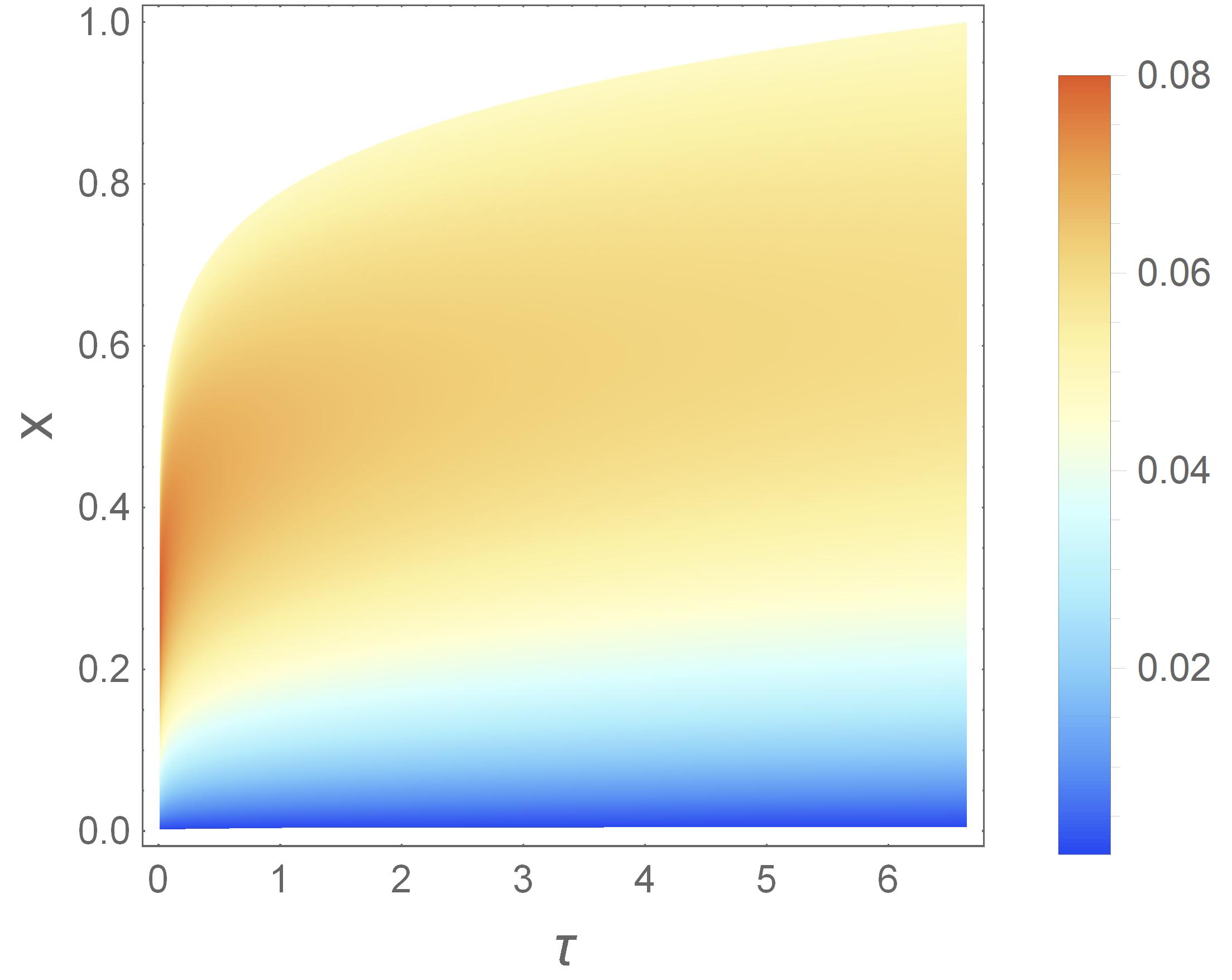}
		\caption{Absolute error generated by the old version of the numerical scheme for 
			$\lambda=\frac{2}{3}$ and $\alpha=0.25$.}
		\label{fig:9}
	\end{center}
\end{figure}

\begin{figure}[h!t]
	\begin{center}
		\includegraphics[scale=0.75]{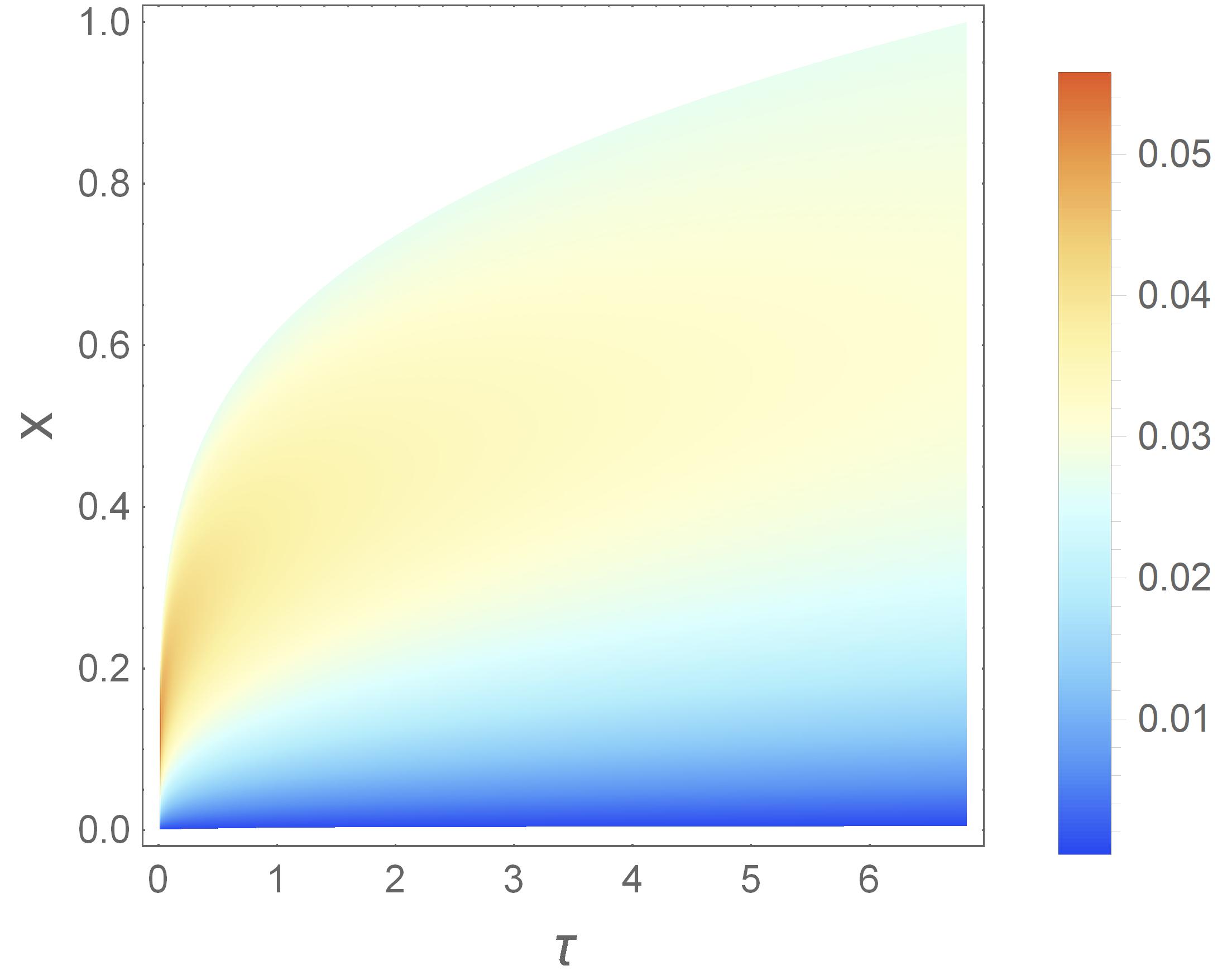}
		\caption{Absolute error generated by the old version of the numerical scheme for 
			$\lambda=\frac{1}{3}$ and $\alpha=0.5$.}
		\label{fig:14}
	\end{center}
\end{figure}

\begin{figure}[h!t]
	\begin{center}
		\includegraphics[scale=0.75]{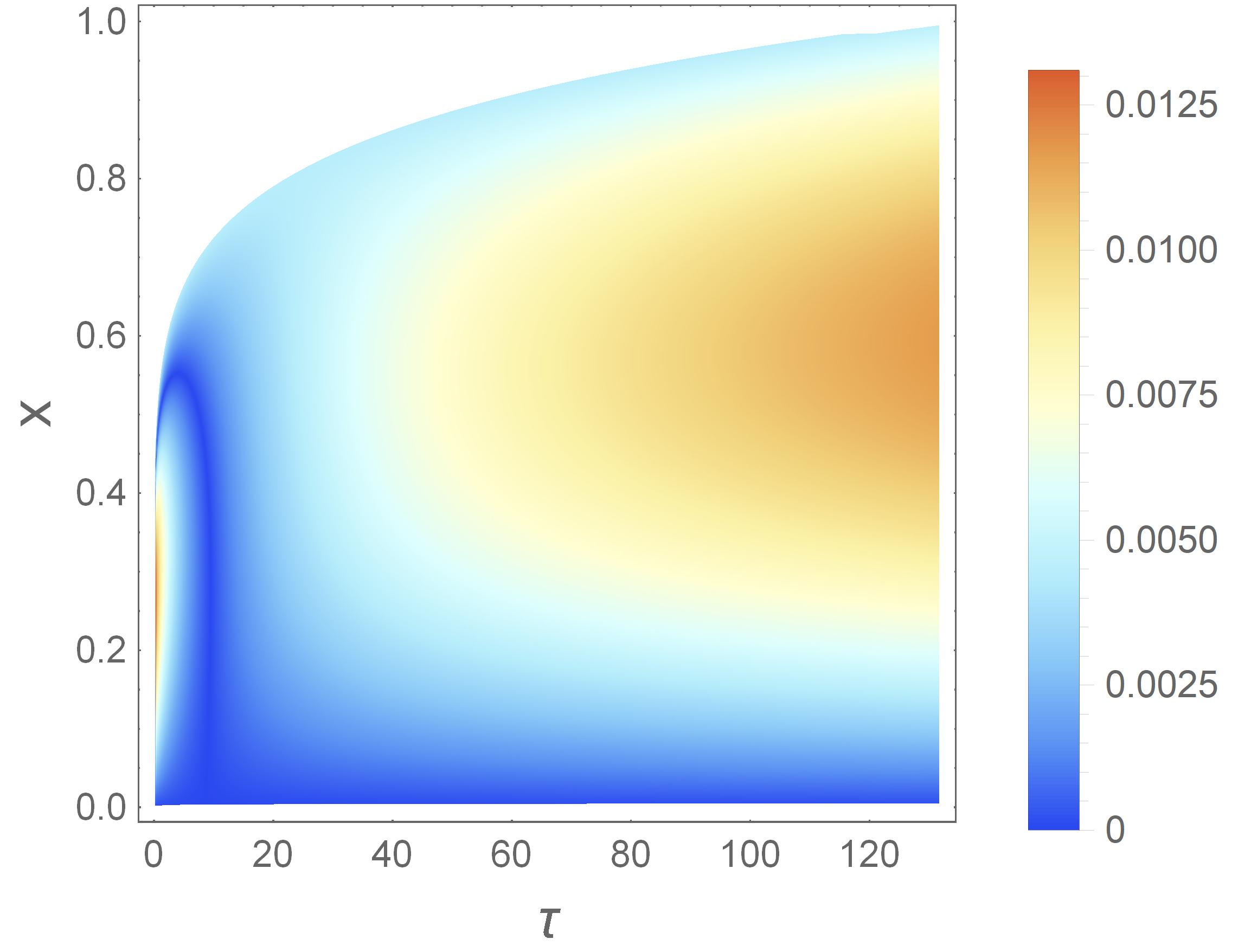}
		\caption{Absolute error generated by the new version of the numerical scheme for 
			$\lambda=\frac{1}{3}$ and $\alpha=0.25$.}
		\label{fig:5}
	\end{center}
\end{figure}

\begin{figure}[h!t]
	\begin{center}
		\includegraphics[scale=0.75]{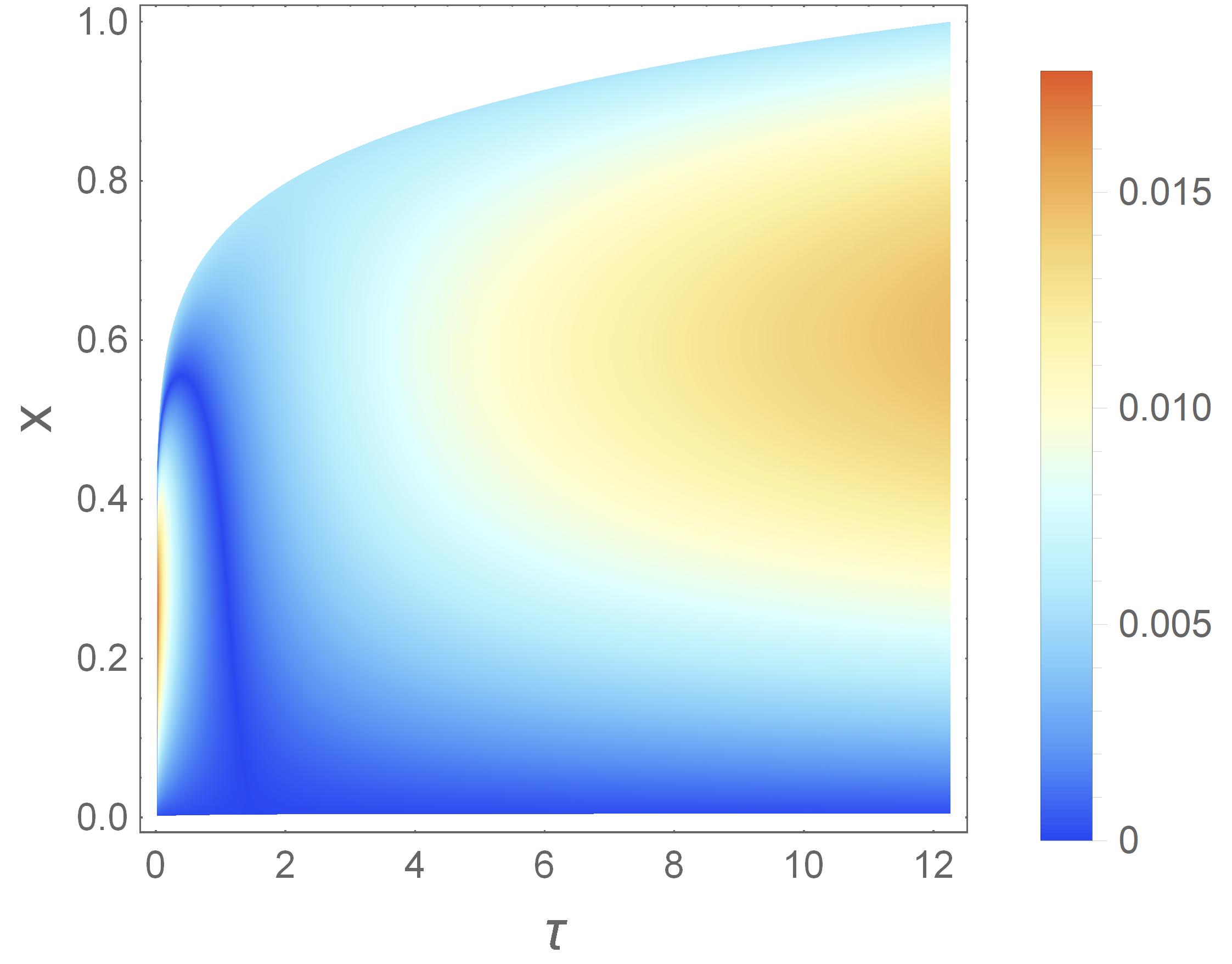}
		\caption{Absolute error generated by the new version of the numerical scheme for 
			$\lambda=\frac{2}{3}$ and $\alpha=0.25$.}
		\label{fig:10}
	\end{center}
\end{figure}

\begin{figure}[h!t]
	\begin{center}
		\includegraphics[scale=0.75]{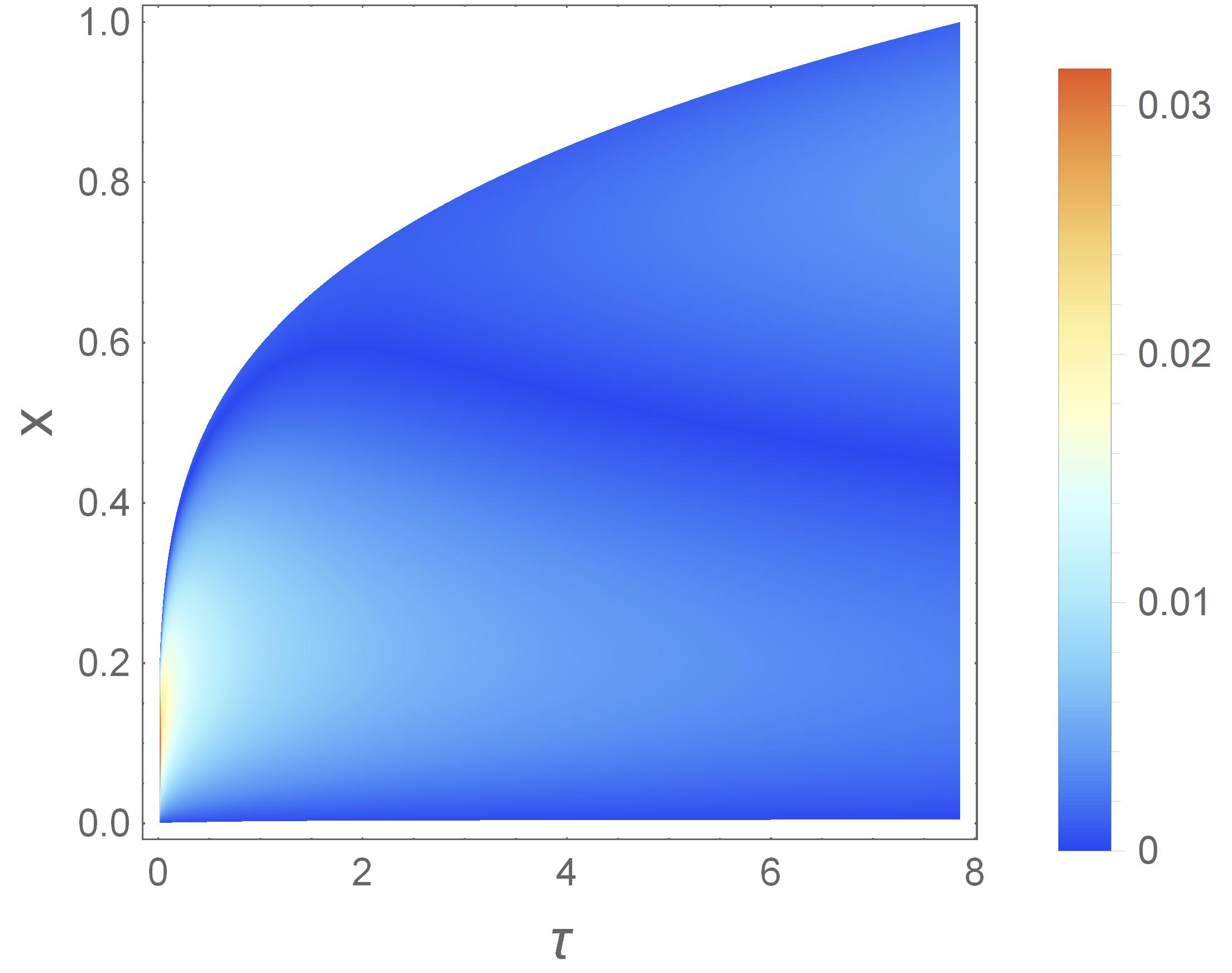}
		\caption{Absolute error generated by the new version of the numerical scheme for 
			$\lambda=\frac{1}{3}$ and $\alpha=0.5$.}
		\label{fig:15}
	\end{center}
\end{figure}

The method proposed in the paper is a modification of the numerical method described in
\cite{Bla15}, therefore we adopted the same values of the model parameters for validation
of the new version of the method. Two variants of the meshes were used in the calculations:
grid with a division into $81\times 241$ and $201\times 801$ nodes. The smaller mesh was used
to generate the results i Tables \ref{tab:p2} and \ref{tab:p3}, while the mesh with higher nodes
density allowed to generate the data illustrated in Figures \ref{fig:1}-\ref{fig:15}. Let
us pay attention to an important fact, in the validation process of the neural network $\phi$,
we assumed sample imputes different from the samples in the training set.

In Table \ref{tab:p} we collected the values of the coefficient $p$ obtained from transcendental 
equation (\ref{eq:trans}) using the Newton algorithm.
\begin{table}[h!t]
	\centering
	\caption{Value of parameter \textit{p} obtained from equation (\ref{eq:trans}).}
	\label{tab:p}
	\begin{tabular}{rcccc}
		\hline
		$\lambda$ & $\alpha=0.25$& $\alpha=0.5$ & $\alpha=0.75$ & $\alpha=1$\\
		\hline
		1/3 & 0.546438 & 0.598238 & 0.669592 & 0.77614 \\
		2/3 & 0.736836 & 0.808016 & 0.90623 & 1.05134 \\
		1   & 0.871649 & 0.956298 & 1.07232 & 1.24014 \\
		\hline
	\end{tabular}
\end{table}

Table \ref{tab:p2} shows the values of the coefficient $p$ obtained by the old version 
of the method. A comparison of the values in Tables \ref{tab:p2} and \ref{tab:p3} with 
Table \ref{tab:p} clearly indicates that the scheme supported by the neural network gives 
much more accurate results.
\begin{table}[h!t]
	\centering
	\caption{Values of parameter $p$ obtained using the old version of the scheme.}
	\label{tab:p2}
	\begin{tabular}{rcccc}
		\hline
		$\lambda$ & $\alpha=0.25$& $\alpha=0.5$ & $\alpha=0.75$ & $\alpha=1$\\
		\hline
		1/3 & 0.584814 & 0.618749 & 0.674536 & 0.774755 \\
		2/3 & 0.78916 & 0.822729 & 0.895849 & 1.050512 \\
		1   & 0.928442 & 0.957251 & 1.041601 & 1.240087 \\
		\hline
	\end{tabular}
\end{table}

\begin{table}[h!t]
	\centering
	\caption{Values of parameter $p$ obtained using the new version of the scheme.}
	\label{tab:p3}
	\begin{tabular}{rcccc}
		\hline
		$\lambda$ & $\alpha=0.25$& $\alpha=0.5$ & $\alpha=0.75$ & $\alpha=1$\\
		\hline
		1/3 & 0.543311 & 0.597387 & 0.660986 & 0.774755 \\
		2/3 & 0.730993 & 0.803931 & 0.892553 & 1.050512 \\
		1   & 0.862646 & 0.948828 & 1.041601 & 1.240087 \\
		\hline
	\end{tabular}
\end{table}
In the further part of the paper, we will show detailed results obtained for 
$\alpha\in\{0.25,0.5\}$ and $\lambda\in\{\frac{1}{3},\frac{2}{3}\}$, because according to
Figure \ref{fig:0} and data collected in the tables above, the new numerical scheme has a 
large advantage over the old one, especially for small values of $\alpha$ and $\lambda$.

In Figure \ref{fig:1}, \ref{fig:6} and \ref{fig:11}, we presented the drug concentration 
profiles obtained by the old version of the front fixing method. In all cases, the graphs were 
plotted for the time when the solid core of the drug was completely dissolved.
\newpage${}$
Similar results were presented in graph \ref{fig:2}, \ref{fig:7} and 
\ref{fig:12}, but were received by the method supported by an artificial neural network.


A boundary surface on which dissolution occurs, moves according to the graphs shown in 
Figure \ref{fig:3}, \ref{fig:8} and \ref{fig:13}. An analytical solution of the function $s$
has been marked with a black line. The green and red lines represent the numerical solution 
obtained with the old version of the front fixing method and the new version respectively. Graph 
analysis in all three cases clearly indicates the advantage of the method supported by an 
artificial neural network.

\newpage
Figure \ref{fig:4}, \ref{fig:9} and \ref{fig:14} shows the distribution of absolute errors
in the domain $[0, s(\tau^*)]\times [0,\tau^*]$ generated by the old version of the front 
fixing method. 

\newpage${}$\newpage
In the case of graphs in Figure \ref{fig:5}, \ref{fig:10} and \ref{fig:15}
we notice that the new method generates much smaller absolute errors. The analysis of 
Figures \ref{fig:4}-\ref{fig:15} also leads to an important observation. In all considered
cases, we note the highest absolute error values for small values of variable $\tau$. 

\newpage
The essence of the front fixing method is the use of a new spatial variable (\ref{eq:sim})
which has singularity at $\tau=0$ and is not defined at this point. Using a numerical approach
consisting in taking as a starting point a very small positive number we generate some error at 
the beginning of calculations. This initial error is transferred to subsequent time layers of the mesh,
indirectly by the preceding time layers, and directly through the first time layer in proportion to a 
certain weight. It seems that the numerical method supported by an artificial neural network allows us
to limit the adverse effect of the results from the initial time layers.
\newpage
\section{Conclusions}
In this paper we proposed a new version of the front fixing method supported by an artificial neural
network. The new scheme allows us to obtain accurate results especially in cases where the previous
version ot the method generated relatively large errors, i.e. for $\alpha<0.5$. It should be noted 
that the neural network presented in the paper is only a simple example of the application of 
artificial intelligence in the optimization process of numerical integration. The use of a more complex
network could hypothetically improve results, but at the expense of simplicity and ease of use.
The presented approach can also be extended to the two-phase fractional Lam\'{e}-Clapeyron-Stefan 
problem.
\bibliography{References}

\begin{thebibliography}{10}

\bibitem{Wee96}
E.R. Weeks, J.S. Urbach, and L.~Swinney.
\newblock Anomalous diffusion in asymmetric random walks with a
  quasi-geostrophic flow example.
\newblock {\em Physica D: Nonlinear Phenomena}, 97:291--310, 1996.

\bibitem{Sol93}
T.H. Solomon, E.R. Weeks, and H.L. Swinney.
\newblock Observations of anomalous diffusion and {L}\'{e}vy flights in a
  2-dimensional rotating flow.
\newblock {\em Physical Review Letters}, 71:3975--3979, 1993.

\bibitem{Hum10}
N.~E. Humphries, N.~Queiroz, J.~R.~M. Dyer, N.~G. Pade, M.~K. Musyl, K.~M.
  Schaefer, D.~W. Fuller, J.~M. Brunnschweiler, T.~K. Doyle, J.~D.~R. Houghton,
  G.~C. Hays, C.~S. Jones, L.~R. Noble, V.~J. Wearmouth, E.~J. Southall, and
  David~W. Sims.
\newblock Environmental context explains {L}\'{e}vy and {B}rownian movement
  patterns of marine predators.
\newblock {\em Nature}, 465:1066--1069, 2010.

\bibitem{Kos05a}
T.~Koszto{\l}owicz, K.~Dworecki, and S.~Mr{\'o}wczy{\'n}ski.
\newblock How to measure subdiffusion parameters.
\newblock {\em Physical Review Letters}, 94:170602, 2005.

\bibitem{Kos05b}
T.~Koszto{\l}owicz, K.~Dworecki, and S.~Mr{\'o}wczy{\'n}ski.
\newblock Measuring subdiffusion parameters.
\newblock {\em Physical Review E}, 71:041105, 2005.

\bibitem{Met00}
R.~Metzler and J.~Klafter.
\newblock The random walk:s guide to anomalous diffusion: a fractional dynamics
  approach.
\newblock {\em Physics Reports}, 339:1--77, 2000.

\bibitem{Met04}
R.~Metzler and J.~Klafter.
\newblock The restaurant at the end of the random walk: recent developments in
  the description of anomalous transport by fractional dynamics.
\newblock {\em Journal of Physics A: Mathematical and General}, 37:161--208,
  2004.

\bibitem{Ste91}
J.~Stefan.
\newblock Uber die theorie der eisbildung, insbesondere uber die eisbildung im
  polarmeere.
\newblock {\em Annalen der Physik und Chemie}, 278(2):269--286, 1891.

\bibitem{Cra84}
J.~Crank.
\newblock {\em Free and Moving Boundary Problems}.
\newblock Clarendon Press, Oxford, 1984.

\bibitem{Hil87}
J.M. Hill.
\newblock {\em One-Dimensional Stefan Problems: An Introduction}.
\newblock Longman Scientific and Technical, New York, 1987.

\bibitem{Rub71}
L.I. Rubinstein.
\newblock {\em The Stefan Problem}.
\newblock American Mathematical Society, 1971.

\bibitem{Gup03}
S.C. Gupta.
\newblock {\em The Classical Stefan Problem. Basic Concepts, Modeling and
  Analysis}.
\newblock Elsevier, Amsterdam, 2003.

\bibitem{Ozi93}
M.~N. Ozisik.
\newblock {\em Heat Conduction}.
\newblock Wiley, 2 edition, 1993.

\bibitem{Liu04}
J.~Liu and M.~Xu.
\newblock An exact solution to the moving boundary problem with fractional
  anomalous diffusion in drug release devices.
\newblock {\em Zeitschrift fur Angewandte Mathematik und Mechanik}, 84:22--28,
  2004.

\bibitem{Vol10}
V.R. Voller.
\newblock An exact solution of a limit case stefan problem governed by a
  fractional diffusion equation.
\newblock {\em International Journal of Heat and Mass Transfer}, 53:5622--5625,
  2010.

\bibitem{Vol13}
V.R. Voller and F.~Falcini.
\newblock Two exact solutions of a stefan problem with varying diffusivity.
\newblock {\em International Journal of Heat and Mass Transfer}, 58:80--85,
  2013.

\bibitem{Sin11}
J.~Singh, P.K. Gupta, and K.N. Rai.
\newblock Homotopy perturbation method to space-time fractional solidification
  in a finite slab.
\newblock {\em Applied Mathematical Modelling}, 35:1937--1945, 2011.

\bibitem{Raj13}
Rajeev and M.S. Kushwaha.
\newblock Homotopy perturbation method for a limit case stefan problem governed
  by fractional diffusion equation.
\newblock {\em Applied Mathematical Modelling}, 37:3589--3599, 2013.

\bibitem{Xic08}
Xicheng Li, Mingyu Xu, and Shaowei Wang.
\newblock Scale-invariant solutions to partial differential equations of
  fractional order with a moving boundary condition.
\newblock {\em Journal of Physics A: Mathematical and Theoretical}, 41:155202,
  2008.

\bibitem{Ros13}
S.~Roscani and E.~Marcus.
\newblock Two equivalent stefan?s problems for the time fractional diffusion
  equation.
\newblock {\em Fractional Calculus and Applied Analysis}, 16(4):802--815, 2013.

\bibitem{Ros16}
S.D. Roscani.
\newblock Hopf lemma for the fractional diffusion operator and its application
  to a fractional free-boundary problem.
\newblock {\em Journal of Mathematical Analysis and Applications},
  434(1):125--135, 2015.

\bibitem{Ala06}
M.~A. Al-Alaoui.
\newblock Al-alaoui operator and the $\alpha$-approximation for discretization
  of analog systems.
\newblock {\em Facta Universitatis Ser.: Elec. Energ.}, 19(1):143--146, 2006.

\bibitem{Liu09}
Liu Junyi and Xu~Mingyu.
\newblock Some exact solutions to stefan problems with fractional differential
  equations.
\newblock {\em Journal of Mathematical Analysis and Applications},
  351:536--542, 2009.

\bibitem{Xia15}
X.~Gao, X.~Jiang, and S.~Chen.
\newblock The numerical method for the moving boundary problem with
  space-fractional derivative in drug release devices.
\newblock {\em Applied Mathematical Modelling}, 39:2385--2391, 2015.

\bibitem{Vol14}
V.R. Voller.
\newblock Fractional stefan problems.
\newblock {\em International Journal of Heat and Mass Transfer}, 74:269--277,
  2014.

\bibitem{Bla15}
M.~B{\l}asik and M.~Klimek.
\newblock Numerical solution of the one phase 1d fractional {S}tefan problem
  using the front fixing method.
\newblock {\em Mathematical Methods in the Applied Sciences},
  38(15):3214--3228, 2015.

\bibitem{Bla14}
M.~B{\l}asik.
\newblock Numerical scheme for one-phase 1d fractional stefan problem using the
  similarity variable technique.
\newblock {\em Journal of Applied Mathematics and Computational Mechanics},
  13(1):13--21, 2014.

\bibitem{Kil06}
A.A. Kilbas, H.M. Srivastava, and J.J. Trujillo.
\newblock {\em Theory and Applications of Fractional Differential Equations}.
\newblock Elsevier, Amsterdam, 2006.

\bibitem{Bar96}
G.I. Barenblatt.
\newblock {\em Scaling, Self-Similarity, and Intermediate Asymptotics}.
\newblock Cambridge University Press, Cambridge, 1996.

\bibitem{Bar03}
G.I. Barenblatt.
\newblock {\em Scaling}.
\newblock Cambridge University Press, Cambridge, 2003.

\bibitem{Ala93}
M.~A. Al-Alaoui.
\newblock Novel digital integrator and differentiator.
\newblock {\em Electronics Letters}, 29(4):376--378, 1993.

\bibitem{Ala08}
M.~A. Al-Alaoui.
\newblock Al-alaoui operator and the new transformation polynomials for
  discretization of analogue systems.
\newblock {\em Electrical Engining}, 90(6):455--467, 2008.

\end{thebibliography}
\bibliographystyle{unsrt}
\end{document}